\documentclass[11pt,reqno]{amsart}

\usepackage[english]{babel}
\usepackage{a4wide}
\usepackage{graphicx}
\usepackage{amsmath,amssymb,amsthm}
\usepackage[ansinew]{inputenc}
\usepackage{listings}
\usepackage{rotating}
\usepackage{hhline}
\usepackage{multirow}
\usepackage{enumerate}
\usepackage[bookmarks]{}
\usepackage{amsfonts}   
\usepackage{tikz}
\usepackage{url}

\newtheorem{theorem}{Theorem}[section]
\newtheorem{lemma}[theorem]{Lemma}

\theoremstyle{definition}
\newtheorem{definition}[theorem]{Definition}

\theoremstyle{remark}

\numberwithin{equation}{section}

\newcommand{\floor}[1]{\left\lfloor #1 \right\rfloor}

\newcommand{\haa}[1]{( #1 )}
\newcommand{\bighaa}[1]{\big( #1 \big)}

\newcommand{\bigaccv}[2]{\big\{ #1 \, \big| \, #2 \big\}}

\newcommand{\Bighaa}[1]{\Big( #1 \Big)}

\newcommand{\biggaccv}[2]{\bigg\{ #1 \, \bigg| \, #2 \bigg\}}

\newcommand{\lrhaa}[1]{\left( #1 \right)}

\newcommand{\lrabs}[1]{\left| #1 \right|}

\newcommand{\pa}[2]{\frac{\partial #1 }{\partial #2 }}

\newcommand{\ga}[2]{\frac{d #1 }{d #2 }}

\newcommand{\R}[1]{\mathbb{R}^{ #1 }}
\newcommand{\PP}{\mathcal{P} (\R{})}
\newcommand{\N}[1]{\mathbb{N}^{ #1 }}

\newcommand{\intabx}[4]{ \int_{ #1 }^{ #2 } { #3 } \hspace{1mm} d { #4 } }
\newcommand{\intabxnohs}[4]{ \int_{ #1 }^{ #2 } { #3 } d { #4 } }

\newcommand{\ddollar}[1]{$\displaystyle{ #1 }$} 

\newcommand{\lrar}{\hspace{1mm} \Leftrightarrow \hspace{1mm}}

\newcommand\Tstrut{\rule{0pt}{4.4ex}}
\newcommand\Bstrut{\rule[-3ex]{0pt}{0pt}}

\DeclareMathOperator{\Dom}{Dom}
\newcommand{\Eintro}{E}
\newcommand{\Enintro}{E_n}
\newcommand{\Enintrotilde}{\tilde{E}_n}
\newcommand{\eps}{\varepsilon}
\newcommand{\evi}{EVI}
\newcommand{\evin}{EVI$_n$}

\newcommand{\Ptwo}{\mathcal{P}_2 (\R{})}

\newcommand{\Veff}{V_{\operatorname{eff}}}
\newcommand{\Wdist}[2]{W_2 \left( #1 , #2 \right) }
\newcommand{\Wdistnobraces}{W_2}
\newcommand{\Wdistq}[2]{W_2^2 \left( #1 , #2 \right) }

\newcommand{\weakto}{\rightharpoonup}

\begin{document}

\title{Upscaling of the dynamics of dislocation walls}


\author{Patrick van Meurs}
\address{Centre for Analysis, Scientific computing and Applications (CASA) \\
  Department of Mathematics and Computer Science \\
  Eindhoven University of Technology \\
  P.O. Box 513 \\
  5600 MB Eindhoven \\ 
  The Netherlands}
\email{p.j.p.v.meurs@tue.nl}
\thanks{The first author is supported by NWO Complexity grant 645.000.012.}

\author{Adrian Muntean}
\address{Centre for Analysis, Scientific computing and Applications (CASA) \\
Institute for Complex Molecular Systems (ICMS)\\
  Department of Mathematics and Computer Science \\
  Eindhoven University of Technology \\
  P.O. Box 513 \\
  5600 MB Eindhoven \\ 
  The Netherlands}





\begin{abstract}
We perform the discrete-to-continuum limit passage for a microscopic model describing the time evolution of dislocations in a one dimensional setting.
This answers the related open question raised by Geers et al. in \cite{GeersPeerlingsPeletierScardia13}. The proof of the upscaling procedure (i.e.~the discrete-to-continuum passage) relies on the gradient flow structure of both the discrete and continuous energies of dislocations set in a suitable evolutionary variational inequality framework. Moreover, the convexity and $\Gamma$-convergence of the respective energies are properties of paramount importance for our arguments.
\end{abstract}

\maketitle

\section{Introduction}
\label{}

\begin{figure}[ht]
\centering
\includegraphics[scale=0.4]{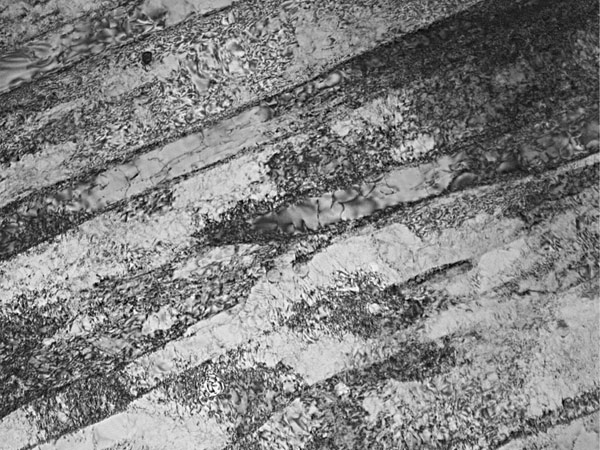}
\caption{Typical dislocation band structures. Courtesy of $\copyright$ Fraunhofer IWS Dresden, Germany \cite{website}.}\label{fig:dlc:pic:bands}
\end{figure}

Plasticity is facilitated at meso- and macro-scales by the movement of line defects (also called dislocations, see Figure \ref{fig:dlc:pic:bands}) in the crystalline structure of the material \cite{vanGoethem, HirthLothe82, HullBacon01}. Describing plasticity at observable  scales (typically $1 - 100$ $\mu$m) is a difficult matter especially due to the cumbersome averaging of the two-scale repulsion interactions (acting at small scales $10 - 100$ $\eta$m) inherently occurring in the discrete dislocation dynamics  to scenarios where dislocations can be well described in terms of a (continuous) density. Even for simplified cases consisting of a two dimensional microscopic model with dislocation lines assumed to be straight and parallel, a general agreement is not yet reached, and hence, conceptually different macroscopic mathematical models exist; see  e.g. \cite{DengElAzab09, ElAzab00, Groma97, GromaBalogh99, GromaCsikorZaiser03, LimkumnerdVan-der-Giessen08, RoyAcharya06, Acharya10}, and \cite{KooimanHutterGeers14}. 

In this paper, we shed light on this issue by performing the rigorous upscaling of the discrete dislocation dynamics for five parameter regimes, rediscovering and classifying this way most of the existing macroscopic models for the evolution of densities of dislocation walls. To capture two-dimensional effects, the basic model of dislocation walls (see Figure~\ref{fig configuration model})
has been considered in \cite{GeersPeerlingsPeletierScardia13, Hall11, RoyPeerlingsGeersKasyanyuk08, VanMeursMunteanPeletier14} as starting point for the discrete-to-continuum limit passage in the static case. The focus here lies exclusively on the dynamic case. 
Note that the discrete-to-continuum passage in the dynamics context has been performed until now only when restricting the dislocations to lie on one slip plane \cite{ForcadelImbertMonneau09}, where the dislocations are assumed to never come closer than a fixed spacing, say $c/n$, so that, in particular, they cannot pile-up as emphasised in the engineering literature (compare e.g. \cite{RoyPeerlingsGeersKasyanyuk08}). In this framework, we relax such restrictions and cover the correct scaling regimes for one-dimensional configurations of dislocations allowing for the piling up of the dislocations at one of the boundaries.

\section{Microscopic dynamics of dislocation walls}
\label{sec dyns of dlc walls}

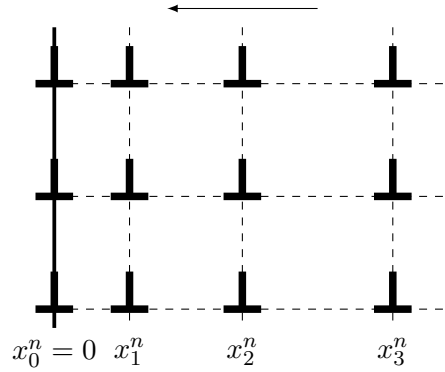
\begin{figure}
\centering
\begin{tikzpicture}[scale=0.5, >= latex]
    \draw[line width = 0.5mm] (0, 7.5) -- (0, -0.5) node [below] {${x}_0^n = 0$};

    \draw[dashed] (2, 7.5) -- (2, -0.5) node [below] {${x}_1^n$};
    \draw[dashed] (5, 7.5) -- (5, -0.5) node [below] {${x}_2^n$};
    \draw[dashed] (9, 7.5) -- (9, -0.5) node [below] {${x}_3^n$};

    \foreach \y in {0,1,2}
        {
        \draw[dashed] (0, 3*\y) -- (10.5, 3*\y);

        \draw[line width = 1mm] (0, 3*\y) -- (0, 3*\y + 1);
        \draw[line width = 1mm] (0 - 0.5, 3*\y) -- (0 + 0.5, 3*\y);
        }
    \foreach \x in {2,5,9}
      \foreach \y in {0,1,2}
        {
        \draw[line width = 1mm] (\x, 3*\y) -- (\x, 3*\y + 1);
        \draw[line width = 1mm] (\x - 0.5, 3*\y) -- (\x + 0.5, 3*\y);
        }

    \draw[<-] (3,8)--(7,8); 
\end{tikzpicture}
\caption{Microscopic geometry showing the vertically periodic configuration of $4$ dislocation walls ($n = 3$). The arrow illustrates the applied constant load.}
\label{fig configuration model}
\end{figure}

We consider dislocations arranged equidistantly in $n+1$ vertical walls, resulting in a vertically periodic setting. Figure \ref{fig configuration model} shows a schematic picture of this configuration. We consider the dynamics of the horizontal positions of these walls (labelled by $x^n_i$) as described by a linear drag law, which states that the velocity of each wall is linearly proportional to the force acting on it. For any wall (other than the most left one) this force equals the sum of an applied load (constant in time) and all the repulsive forces exerted by all other walls. The left-most wall at $x_0^n = 0$ is fixed in time, which models an impenetrable barrier.

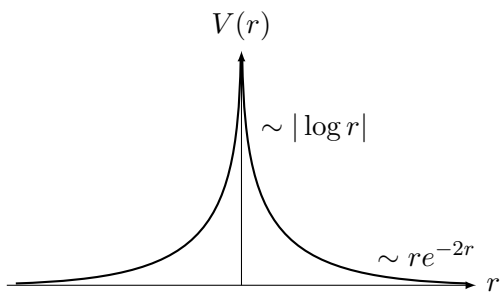
\begin{figure}
\begin{tikzpicture}[scale=0.6, >= latex]
\draw[->] (-5.2,0) -- (5.2,0) node[right] {$r$};
\draw[->] (0,0) -- (0,5.2) node[above] {$V(r)$};
\draw[thick] (0.018,5) .. controls (0.09,1) and (1,0.176) .. (5,0.04);
\draw[thick] (-0.018,5) .. controls (-0.09,1) and (-1,0.176) .. (-5,0.04);
\draw (0.2, 3.5) node[right] {$\sim |\log r|$};
\draw (4.1, 0.7) node {$\sim r e^{-2r}$};
\end{tikzpicture}
\caption{The shape of the interaction potential $V$. Note the two characteristic length scales of repulsion.}
\label{fig:V}
\end{figure}

The dimensionless form of the microscopic model is as follows: 
\begin{equation} \label{for:disc:evo:eqn:in:Rn}
  \left\{ \begin{aligned}
    \displaystyle \ga{}t x^n (t) &= - n \nabla \Enintro \lrhaa{ x^n (t) }, \quad t > 0 \\
    \displaystyle x^n (0) &= x_0^n, \\
  \end{aligned} \right.
\end{equation}
where
\begin{equation}
	x^n = \lrhaa{x^n_i}_{i=1}^n \in \Omega_n := \bigaccv{ (x_1, \ldots, x_n) \in (0, \infty)^n }{ x_1 < x_2 < \ldots < x_n },
\end{equation}
and
\begin{gather} \label{for:defn:En}
  \Enintro (x^n) 
  := 
  \frac{\alpha_n}n \sum_{k=1}^n \sum_{j = 0}^{n - k} V \bighaa{ n \alpha_n \bighaa{ x_{j+k}^n - x_j^n } } 
  + \frac1n \sum_{i = 1}^{n} x_i^n,
\end{gather}
in which $n$ and $\alpha_n$ are the only two dimensionless parameters. The potential $V : \R{} \rightarrow \R{}$ is given by
\begin{equation} \label{for defn V}
	V(r) := r \coth r - \log \lrabs{ \sinh r } - \log 2
\end{equation}
and illustrated in Figure \ref{fig:V}. We refer the reader to \cite{GeersPeerlingsPeletierScardia13, VanMeursMunteanPeletier14} for a derivation of \eqref{for:defn:En}. The parameter $\alpha_n$ is proportional to the \emph{aspect ratio} between the typical horizontal distance and the vertical distance between two neighbouring dislocations in Figure \ref{fig configuration model}. Due to the two characteristic length scales of $V$, the asymptotic behaviour of $\alpha_n$ influences the behaviour of the interactions for large $n$. This effect is central in \cite{GeersPeerlingsPeletierScardia13}. In fact, when $\alpha_n \ll 1/n$ or $\alpha_n \gg 1$, the scaling of $E_n$ is different from \eqref{for:defn:En}. We refer to Appendix \ref{app:precise:scaling} for the details.

Since $\Enintro \in C^\infty (\Omega_n)$ is bounded from below, strictly convex, and has compact level sets in $\Omega_n$, it holds that \eqref{for:disc:evo:eqn:in:Rn} has a unique solution $x^n \in C^1 ([0, \infty); \Omega_n)$ for any initial value in $\Omega_n$. 

Before addressing the upscaling of the gradient flow \eqref{for:disc:evo:eqn:in:Rn}, it is convenient to cast the problem into a form (i) without a gradient, and (ii) with a structure given in an $n$-independent underlying space. To accomplish (i) we take any $y \in \R n$, multiply both sides of \eqref{for:disc:evo:eqn:in:Rn} by $(x^n - y)$, and exploit the convexity of $E_n$ to conclude that
\begin{equation} \label{for:disc:evo:eqn:in:Rn:weak}
  \frac12 \ga{}t \| x^n - y \|^2 
  = 
  - n (x^n - y) \nabla \Enintro ( x^n ) 
  \leq 
  n \bighaa{ \Enintro (y) - \Enintro (x^n) }, 
  \quad
  \textup{for all } y \in \R n .
\end{equation}
To satisfy (ii), we interpret the inequality \eqref{for:disc:evo:eqn:in:Rn:weak} in the space of probability measures with finite second moments, denoted by $\Ptwo$. To this aim we define $\Enintrotilde : \Ptwo \rightarrow (-\infty, \infty]$ by 
\begin{equation}
	\Enintrotilde (\mu^n) = \left\{ \begin{array}{ll}
	  \displaystyle \Enintro (x), & \mu^n \in \Dom \Enintrotilde, \\
	  \displaystyle \infty, & \textup{otherwise,} \\
	\end{array} \right.
\end{equation}
where
\begin{equation} \label{for:defn:dom:En:tilde}
	\Dom \Enintrotilde := \biggaccv{ \mu \in \Ptwo } { \exists \; x \in \Omega_n : \mu^n = \frac 1n \sum_{i=1}^n \delta_{x_i} }.
\end{equation}
We remove the tildes for ease of notation. Then it immediately follows from \eqref{for:disc:evo:eqn:in:Rn:weak} that
\begin{equation} \label{for:evin}
  \left\{ \begin{aligned}
  	  \displaystyle \frac 12 \ga{}t \Wdistq{\mu_t^n}{\nu} &\leq \Enintro (\nu) - \Enintro (\mu_t^n), \quad \textup{for all } t>0, \: \nu \in \Dom \Enintro , \\
  	  \displaystyle \lim_{t \downarrow 0} \Wdist{\mu_t^n}{\mu_0^n} &= 0,  \\
  	\end{aligned} \right.	
\end{equation}
where $\Wdistnobraces$ is the Wasserstein-$2$ metric (see e.g. \cite{AmbrosioGigliSavare08} for its definition and basic properties). Following the terminology from \cite{AmbrosioGigliSavare08}, the inequality in \eqref{for:evin} is called discrete \emph{evolution variational inequality} , say (\evin). Theorem 4.0.4.iii from \cite{AmbrosioGigliSavare08} states the desired well-posedness (see Lemma \ref{lem:AGS} below for the details). The application of such  powerful result requires the underlying metric space to be complete. This is why we are considering $(\Ptwo, \Wdistnobraces)$ rather than, for example, the usual space of probability measures $\PP$ (without any restriction on any of their moments) equipped with the narrow topology.  

\section{Macroscopic dynamics of the dislocation wall density}
\label{sec dyns of dlc wall dens}

With the \evin{} (see \eqref{for:evin}) in hand, we address the question raised in the introduction: ``What is the limit evolution equation as $n \rightarrow \infty$?". An obvious candidate for the energy structure of this limit  evolution equation is given in \cite{GeersPeerlingsPeletierScardia13}, where it is proven that $\Enintro$ $\Gamma$-converges to an explicitly given functional $\Eintro$ (see Table \ref{tab:E}). The expression of $\Eintro$ depends on the asymptotic behaviour of the parameter $\alpha_n$ with respect to $n \to \infty$. Since our analysis below is independent on the asymptotic behaviour of $\alpha_n$ and the expression for $\Eintro$ which it induces, we do not include it in our notation.

The $\Gamma$-convergence of $\Enintro$ to $\Eintro$ is proven on the space of probability measures $\PP$ equipped with the narrow topology. This is a weak topology defined by testing with continuous and bounded functions. More precisely, $(\mu_n) \subset \PP$ converges narrowly to $\mu$ if
\begin{equation*}
\forall \: \varphi \in C_b (0, \infty) : \intabx0\infty{ \varphi }{\mu_n} \rightarrow \intabx0\infty{ \varphi }{\mu}
\end{equation*}

Due to such a $\Gamma$-convergence result, we expect that the following \emph{evolution variational inequality} (\evi) arises as the natural limit of the EVI$_n$ \eqref{for:evin} when $n$ tends to infinity, viz.
\begin{equation} \label{for:evi}
  \left\{ \begin{aligned}
  	  \frac 12 \ga{}t \Wdistq{\mu_t}{\nu} &\leq \Eintro (\nu) - \Eintro (\mu_t), \quad \textup{for all } t>0, \; \nu \in \Dom \Eintro , \\
  	  \displaystyle \lim_{t \downarrow 0} \Wdist{\mu_t}{\mu_0} &= 0,  \\
  	\end{aligned} \right.
\end{equation}
where $\mu_0 \in \Dom \Eintro$. To show the well-posedness of \eqref{for:evi} we use again Theorem~4.0.4.iii in \cite{AmbrosioGigliSavare08} (see Lemma \ref{lem:AGS} for the application of the general result to our precise setting).

\begin{table}[h!]
\centering
\caption{Expressions for the continuous energy $E$ for the five different scaling regimes of $\alpha_n$. We denote $a := \int_0^\infty V$ and $\Veff (r) := \sum_{k = 1}^\infty V \haa{ k r }$.}
\label{tab:E}
\begin{tabular}{|c|l|}
  \hline
  Regime & Structure of the continuous energy $E$ \\\hline
  \Bstrut \Tstrut $\alpha_n \ll \dfrac1n$ & \ddollar{ E^{1} (\mu) = \intabx 0\infty{x}\mu + \frac{1}{2} \intabxnohs{0}{\infty}{ \intabx{0}{\infty}{ \log \frac1{|x - y|} }{\mu (y)} }{ \mu(x) } } \\
  \Bstrut $n \alpha_n \rightarrow {c}$ & \ddollar{ E^{2} (\mu) = \intabx 0\infty{x}\mu + \frac{c}{2} \intabxnohs{0}{\infty}{ \intabx{0}{\infty}{ V ( {c}(x - y) ) }{\mu (y)} }{ \mu(x) } } \\
  \Bstrut $\dfrac1n \ll \alpha_n \ll 1$ & \ddollar{ E^{3} (\mu) = \intabx 0\infty{x}\mu + \left\{
                     \begin{array}{ll}
                       \displaystyle a \int_0^\infty \rho^2, & \hbox{if $d\mu(x) = \rho(x) dx$,} \\
                       \infty, & \hbox{otherwise}
                     \end{array}
                   \right. } \\
  \Bstrut $\alpha_n \rightarrow c$ & \ddollar{ E^{4} (\mu) = \intabx 0\infty{x}\mu + \left\{
                     \begin{array}{ll}
                       \displaystyle c \intabx{0}{\infty}{ \Veff \Bighaa{ \frac{c}{ \rho(x) } } \rho(x) }{x}, & \hbox{if $d\mu(x) = \rho(x) dx$,} \\
                       \infty, & \hbox{otherwise}
                     \end{array}
                   \right. } \\
  $1 \ll \alpha_n$ & \ddollar{ E^{5} (\mu) = \intabx 0\infty{x}\mu + \left\{
                     \begin{array}{ll}
                      0, & \hbox{if \ddollar{ \ga{\mu}{\mathcal{L}} \leq 1} $\mathcal{L}$-a.e.,} \\
                       \infty, & \hbox{otherwise}
                     \end{array}
                   \right. } \\
  \hline
\end{tabular}
\end{table}

\section{Upscaling procedure. Main result}
\label{sec result}

Our main result states that the evolution of the dislocation wall density that we are looking for is indeed given by \eqref{for:evi}. Moreover, it specifies the topology in which the dislocation wall density dynamics is approximated by the discrete dislocation walls dynamics.

\begin{theorem} (Evolutionary convergence). \label{thm}
Let $T > 0$, $n \in \N{}$, and $(\alpha_n)$ such that either $\alpha_n \ll 1/n, \, \alpha_n \sim 1/n, \, 1/n \ll \alpha_n \ll 1, \, \alpha_n \sim 1$, or $\alpha_n \gg 1/n$. Let $E$ be the $\Gamma$-limit of $E_n$, and $\mu_0 \in \overline{\Dom \Eintro}$. Then for any $\mu^n_0 \in \Dom \Enintro$ such that $\Wdist{\mu_0^n}{\mu_0} \rightarrow 0$, it holds for the solutions $\mu_t^n$ and $\mu_t$ to \eqref{for:evin} and \eqref{for:evi} respectively that $\Wdist{\mu_t^n}{\mu_t} \rightarrow 0$ point-wise for all $t > 0$.

\end{theorem}

\subsection{Proof of Theorem \ref{thm}}

The proof heavily relies on the following two technical results: Lemma \ref{lem:AGS} and Lemma \ref{lem:stability}, both adapted versions of Theorem 4.0.4.iii in \cite{AmbrosioGigliSavare08} and Theorem 6.1 in \cite{AmbrosioSavareZambotti09}\footnote{As noted on page 534 of  \cite{AmbrosioSavareZambotti09}, the proof of their Theorem 6.1 implies a stronger statement which covers Lemma \ref{lem:stability}} respectively. 

\begin{lemma} (Well-posedness of gradient flows). \label{lem:AGS}
Let $\phi : \Ptwo \rightarrow (-\infty, \infty]$ be proper, l.s.c., bounded from below, and geodesically convex. Let $u_0 \in \overline{ \Dom \phi }$. Then the EVI (see \eqref{for:evi}) has exactly one solution among all locally absolutely continuous curves in $\Ptwo$.

\end{lemma}

\begin{lemma} (Stability of gradient flows). \label{lem:stability}
Let $\phi_n, \phi : \Ptwo \rightarrow (-\infty, \infty]$ as in Lemma \ref{lem:AGS}. For the initial data, let $(u_0^n) \subset { \Dom \phi_n }$ be a $\Wdistnobraces$-converging sequence with limit $u_0 \in \overline{ \Dom \phi }$. If $\phi_n$ Mosco-converges to $\phi$ (see \eqref{for:Mosco:conv}), then it holds for the solutions $u_t^n$ and $u_t$ to the related gradient flows that $\Wdist{ u_t^n}{ u_t } \rightarrow 0$ point-wise in $t \in [0, \infty)$.

\end{lemma}

We see that Lemma \ref{lem:AGS} applies to both the discrete and continuous energies $\Enintro$ and $\Eintro$ ($\Eintro$ is strictly convex, because $\Gamma$-convergence conserves strictly convexity). Hence, we are left to show the Mosco-convergence of the energies:
\begin{subequations}
\label{for:Mosco:conv}
\begin{alignat}2
\label{for:Mosco:conv:limf}
&\textup{for all } \mu_n \weakto \mu \textup{ narrowly}, & \liminf_{n\to\infty}  \Enintro (\mu_n) &\geq \Eintro (\mu), \\
&\textup{for all } \mu \textup{ there exists }\mu_n \xrightarrow \Wdistnobraces \mu \textup{ such that}\quad & \limsup_{n\to\infty} \Enintro (\mu_n) &\leq \Eintro (\mu),
\label{for:Mosco:conv:limp}
\end{alignat}
\end{subequations}
The liminf inequality \eqref{for:Mosco:conv:limf} and the limsup inequality \eqref{for:Mosco:conv:limp} are given by Theorem~1 in \cite{GeersPeerlingsPeletierScardia13} with respect to narrow topology. The main idea for strengthening the limsup inequality to hold in the Wasserstein-$2$ metric is to consider $\mu_n = \frac 1n \sum_{i=1}^n \delta_{x_i^n}$ defined by
\begin{equation}\label{eqn:rec:seq}
x_i^n := \inf\left\{ x: \mu([0,x]) \geq \frac i{n+1}\right\},
\end{equation}
as the recovery sequences. The difference with the recovery sequence in \cite{GeersPeerlingsPeletierScardia13} is that the number $i/n$ has to be  substituted by $i/(n+1)$. Proving that this strategy applies, requires some technicalities and computations, which are given in Appendix~\ref{app:limsup}.

\section{Rate of convergence}

Theorem \ref{thm} states the convergence of $\Wdist{\mu_t^n}{\mu_t}$ to $0$. In this section, we illustrate numerically the expected convergence  rate.

We denote by $x^{n_k} (t)$ our numerical approximation\footnote{We have used the `ode15s' solver \cite{ShampineReichelt97} in MATLAB, and defined $x^{n_k} (t)$ by linear interpolation in $t$.} of the solution to \eqref{for:disc:evo:eqn:in:Rn} with $\alpha_{n_k} = 1/\sqrt{n_k}$, 
\begin{equation*}
n_k := \floor{ 10 \cdot 2^{k/2} },
\quad
k = 0,1,\ldots,10,
\quad
\lrhaa{ x_0^{n_k} }_i := 2 \sqrt a \frac i{n_k}, 
\quad
a := \int_0^\infty V.
\end{equation*}
Our choice for the initial guess is motivated by the facts that (i) the dislocation walls tend to place themselves equidistantly in the absence of an external force, and (ii) the support of the minimiser of the continuous energy (Table \ref{tab:E}) equals $[0, 2 \sqrt a]$. We choose $\alpha_{n_k} = 1/\sqrt{n_k}$ such that it corresponds to the physical parameter regime from \cite{Dogge14, EversBrekelmansGeers04}.

Let $\mu^{n_k} (t)$ be the related empirical measure according to \eqref{for:defn:dom:En:tilde}. We are interested in how fast $\gamma_{n_k} := \Wdist{\mu_t^{n_{k+1}}}{\mu_t^{n_k}}$ decreases in $n_k$. With the ansatz $\gamma_{n_k} \sim C n_k^{-p}$, we can estimate $p$ by computing
\begin{equation} \label{for:pnk}
   p_{n_k} := - \frac2{ \log 2 } \log \frac{ \gamma_{n_{k+1}} }{ \gamma_{n_k} }.
\end{equation} 
Table \ref{tab:pnk:sqrt} shows the results for $t = 2^{-6}, 2^{-2}, 2^{2}, \infty$. By ``$t = \infty$" we refer to the numerical approximation $\mu^n_\ast$ of the equilibrium profile, which we obtain by putting the right hand side of \eqref{for:disc:evo:eqn:in:Rn} equal to zero and solving with Newton's method. As the typical time for the gradient flow to reach the equilibrium state is between $1$ and $10$ time units, we have chosen the other time instances such that we see what happens: (i) quickly after the initial configuration, (ii) somewhere halfway in reaching equilibrium, (iii) close to the equilibrium, and (iv) at equilibrium.  

\begin{table}[h!]                    
\centering    
\caption{Values of $p_{n_k}$ obtained from the approximation \eqref{for:pnk} for $\alpha_n = 1/\sqrt n$ at different time instances $t$. The decrease of $p_{n_k}$ in time and $k$ suggests that a boundary layer is being formed at the left of the pile-up of dislocation walls.}             
\label{tab:pnk:sqrt}                          
\begin{tabular}{|c|cccc|}       
\hline                               
 \rule{0pt}{2.4ex} & $t=2^{-6}$ & $t=2^{-2}$ & $t=2^2$ & $t=\infty$ \\                                            
\hline                               
 $n_1$ & 0.936 & 0.925 & 0.755 & 0.758 \\  
                              
 $n_2$ & 1.036 & 1.000 & 0.827 & 0.827 \\  
                              
 $n_3$ & 0.939 & 0.920 & 0.730 & 0.724 \\  
                               
 $n_4$ & 1.040 & 0.993 & 0.800 & 0.806 \\  
                               
 $n_5$ & 0.944 & 0.918 & 0.664 & 0.657 \\  
                               
 $n_6$ & 1.005 & 0.962 & 0.710 & 0.728 \\  
                               
 $n_7$ & 0.985 & 0.940 & 0.664 & 0.650 \\  
                               
 $n_8$ & 0.991 & 0.941 & 0.644 & 0.638 \\  
                               
 $n_9$ & 0.981 & 0.916 & 0.572 & 0.583 \\  
\hline                               
\end{tabular}                                  
\end{table} 

The first column suggests that $p = 1$. This is what we would expect, because  for the initial conditions $\mu^{n_k} (0)$ it is easy to calculate that $p = 1$ exactly. The second and the third column show that $p$ decreases both in time and in the number of dislocation walls, which means that the rate of convergence decreases as well. We expect that this decrease is due to the presence of a boundary layer close to the barrier, i.e.~at $x = 0$. Indeed, column $4$ shows that, for the equilibrium profiles (obtained by minimizing $E_n$), the values of $p_{n_k}$ drop just as fast as $n_k$ increases. This suggests that the decrease in $p$ over time is due to the formation of the boundary layer rather than $\mu_n (t)$ becoming a worse approximation of $\mu (t)$, i.e.~the solution to the continuous gradient flow.

The mathematical analysis of the structure of the boundary layer is still in progress \cite{GarronivanMeursPeletierScardia14prep, Hall11}. We expect that this layer consists of $\mathcal O (1/\alpha_n)$ particles; compare Table~\ref{tab:pnk:sqrt} where the values of $p_{n_k}$ tend to $1/2$ as $n_k$ increases for large $t$. 

Now, we restrict our attention to testing whether this decrease in the convergence rate is indeed due to the boundary layer. To this aim, we consider $\alpha_n = 1/n$, because in this parameter regime the boundary layer is expected to consist of a constant number of dislocation walls. This parameter choice is physically relevant; it corresponds to the transition regime between the classical models in \cite{EshelbyFrankNabarro51, HeadLouat55} and the one proposed in \cite{EversBrekelmansGeers04}.

We perform the same numerical investigations as above, but now with $\alpha_n = 1/n$. The results are shown in Table \ref{tab:pnk:one}. The values of $p_{n_k}$ related to the equilibrium profile (fourth column) are all close to $1$, independent of $n_k$. This implies that there is not a boundary layer consisting of a growing number of dislocation walls, as expected. Indeed, we see from the first three columns that $p_{n_k}$ hardly decreases as time elapses or with respect to  $n_k$. Furthermore, its value is close to $1$, which is the highest possible value for the convergence rate.

\begin{table}[h!]                  
\centering    
\caption{Values of $p_{n_k}$ obtained from the approximation \eqref{for:pnk} for $\alpha_n = 1/n$ at different time instances $t$. All these values being close to $1$ shows that the rate of convergence is close to $1/n$.}             
\label{tab:pnk:one}                          
\begin{tabular}{|c|cccc|}         
\hline                               
 \rule{0pt}{2.4ex} & $t=2^{-6}$ & $t=2^{-2}$ & $t=2^2$ & $t=\infty$ \\                      
\hline                               
 $n_1$ & 0.945 & 0.949 & 0.846 & 0.855 \\  
                              
 $n_2$ & 1.046 & 1.038 & 0.924 & 0.929 \\  
                              
 $n_3$ & 0.946 & 0.940 & 0.881 & 0.882 \\  
                              
 $n_4$ & 1.048 & 1.027 & 0.946 & 0.955 \\  
                               
 $n_5$ & 0.951 & 0.940 & 0.893 & 0.902 \\  
                             
 $n_6$ & 1.018 & 1.001 & 0.953 & 0.950 \\  
                              
 $n_7$ & 0.998 & 0.982 & 0.942 & 0.948 \\  
                               
 $n_8$ & 1.004 & 0.994 & 0.962 & 0.957 \\  
                             
 $n_9$ & 0.995 & 0.982 & 0.948 & 0.960 \\  
\hline                               
\end{tabular}                                
\end{table}

\section{Discussion}
\label{sec:discussion}

For practical applications of the main result stated in Theorem \ref{thm},  it is desired to describe the EVI \eqref{for:evi} in terms of an explicit PDE for any energy $E$ in Table \ref{tab:E}. Next we formally show how to obtain these PDEs. 

In Example 11.2.7 in \cite{AmbrosioGigliSavare08} explicit PDEs are given for the gradient flows corresponding to energy functionals $\mathcal E : \Ptwo \to \R{}$ of the form
\begin{equation} \label{for:E:abstract}
    \mathcal E (\mu)
    :=
    \left\{ \begin{array}{ll}
      \displaystyle 
        \int \Phi (\rho) 
    	+ \int \mathcal F \rho
    	+ \frac12 \int \rho (\mathcal V \ast \rho)
      & \hbox{if $d\mu(x) = \rho(x) dx$,} \\
      \displaystyle \infty & \hbox{otherwise,} \\
    \end{array} \right.
\end{equation}
in which $\Phi, \mathcal F, \mathcal V$ should satisfy a list of assumptions; see Section 10.4.7 in \cite{AmbrosioGigliSavare08} for the details. These gradient flows are given by
\begin{equation} \label{for:PDE:abstract}
  \pa{}{t} \rho_t 
  = 
  \bighaa{ 
    (L_\Phi \circ \rho_t)' 
    + \rho_t \mathcal F' 
    + \rho_t (\mathcal V' \ast \rho_t) 
    }' 
  \qquad
  \text{in $\mathcal D' (\R{} \times (0, \infty))$},
\end{equation}
where $L_\Phi (z) := z \Phi'(z) - \Phi (z)$, together with the conditions of mass conservation of $\rho_t$, integrability of $L_\Phi \circ \rho_t \in L^1_{\text{loc}} (\R{} \times (0, \infty))$, $\rho_t$ having finite second moments, and $\Wdist{\rho_t \mathcal L}{\mu_0} \to 0$ as $t \searrow 0$. The unique solution to the gradient flow of $\mathcal E$ (i.e.~the solution to \eqref{for:evi}) is equal to the unique solution to the PDE \eqref{for:PDE:abstract} with the conditions stated below it (see Theorem~11.2.8 in \cite{AmbrosioGigliSavare08}).

 Our expressions $E^{i}$ as listed in Table \ref{tab:E} are almost of the type \eqref{for:E:abstract}. To see this, we take for any $i$
\begin{equation*}
  \mathcal F (x) 
  :=
  \left\{ \begin{array}{ll}
    \displaystyle \infty & x < 0 \\
    \displaystyle x & x \geq 0 \\
  \end{array} \right.
\end{equation*}
to account for the externally applied force and the condition that the dislocation walls cannot penetrate the barrier at $x = 0$. For $i = 1,2$, we take $\Phi \equiv 0$ and $- \log |\cdot|$ and $c V (c \, \cdot)$ respectively for $\mathcal V$. This choice for $\mathcal V$ does not fit into \eqref{for:E:abstract}, because $\mathcal V$ is assumed to be convex and differentiable, while our choices are not at $x = 0$. For $i = 3,4$, we take $\mathcal V \equiv 0$ and respectively $a z^2$ and $cz\Veff (c/z)$ for $\Phi(z)$. Both these choices satisfy all requirements on $\mathcal V$ and $\Phi$. Case $i = 5$ is too degenerate for \eqref{for:PDE:abstract} to make sense. It is not clear how to make sense of the related flux. As this parameter regime is not of physical interest, we do not consider it further. 

With these choices for $\Phi, \mathcal F, \mathcal V$ we find for $i = 1, \ldots, 4$ respectively the following effective PDEs
\begin{align}
  \pa{}t \rho_t 
  &= 
  \bighaa{ \rho_t + \rho_t \haa{ -\log | \cdot | * \rho_t }' }' 
  && 
  \textup{if } \alpha_n \ll \frac1n, 
  \label{for:evo:eq:1}
  \\
  \pa{}t \rho_t 
  &= 
  \bighaa{ \rho_t + \rho_t \haa{ c V( c \, \cdot ) * \rho_t }' }' 
  && 
  \textup{if } n \alpha_n \rightarrow {c}, 
  \label{for:evo:eq:2}
  \\
  \pa{}t \rho_t 
  &= 
  \lrhaa{ \rho_t + 2 a \rho_t \rho_t' }' 
  && 
  \textup{if } \frac1n \ll \alpha_n \ll 1, 
  \label{for:evo:eq:3}
  \\
  \pa{}t \rho_t 
  &= 
  \lrhaa{ \rho_t + \frac{c^2}{\rho_t^2} \Veff'' \Bighaa{ \frac c{\rho_t} } \rho_t' }' 
  && 
  \textup{if } \alpha_n \rightarrow c, 
  \label{for:evo:eq:4}
\end{align}
with the initial condition satisfying $\Wdist{\rho_t \mathcal L}{\mu_0} \to 0$. The boundary conditions are obtained from the conditions ensuring the conservation of mass. This yields zero flux at $x = 0$, and $\rho_t (x) \to 0$ as $x \to \infty$.

Although Theorem \ref{thm} is restricted to dislocation wall dynamics on a half-line, the proof allows for generalisations to include for instance the case of dislocations moving in finite domains cf.~the setting discussed  in \cite{VanMeursMunteanPeletier14}. The related PDEs are similar to \eqref{for:evo:eq:1} -- \eqref{for:evo:eq:4} equipped with another zero flux boundary condition at the other end of the finite domain.

\section*{Acknowledgements}
We like to express our gratitude towards Mark~Peletier (Eindhoven) for many fruitful discussions. Our thanks also go to R.~Rossi (Brescia) and G.~Savar\'e (Pavia) for providing us with useful advice.

\appendix

\section{Precise scaling for $E_n$}
\label{app:precise:scaling}

As mentioned in Section \ref{sec dyns of dlc walls}, the scaling of the discrete energy $E_n$ is slightly different when $\alpha_n \ll 1/n$ or $\alpha_n \gg 1$. The precise scaling for these two cases is given by
\begin{equation*} 
  \Enintro (x^n) 
  := 
  \frac 1{n^2} \sum_{k=1}^n \sum_{j = 0}^{n - k} V \bighaa{ n \alpha_n \bighaa{ x_{j+k}^n - x_j^n } } 
  - \frac12 \log \frac{e}{ 2 n \alpha_n }
  + \frac1n \sum_{i = 1}^{n} x_i^n,
\end{equation*}
when $\alpha_n \ll 1/n$, and  
\begin{equation*} 
  \Enintro (x^n) 
  := 
  \frac{ \exp \bighaa{2 \haa{ \alpha_n - 1 } } }{ n \alpha_n } \sum_{k=1}^n \sum_{j = 0}^{n - k} V \bighaa{ n \alpha_n \bighaa{ x_{j+k}^n - x_j^n } } 
  + \frac1n \sum_{i = 1}^{n} x_i^n,
\end{equation*}
when $\alpha_n \gg 1$.

\section{Proof of the limsup inequality \eqref{for:Mosco:conv:limp} }
\label{app:limsup}

In tis appendix, we show the details of the proof of the limsup inequality \eqref{for:Mosco:conv:limp}. To do so, we modify the proof of Theorem 1 in \cite{GeersPeerlingsPeletierScardia13} for the first four parameter regimes.  For the regime $\alpha_n \gg 1$ we proceed differently -- here the arguments from Theorem 1.1 in \cite{VanMeursMunteanPeletier14} need to be adapted. To extend these proofs properly to deal with the current scenario, we require some preliminaries to be able to relate the narrow topology to the Wasserstein-$2$ metric:

\begin{definition} (Uniformly integrable second moments).
A set $A \subset \PP$ has \emph{uniformly integrable second moments} if for all $\eps > 0$ there exists an $M > 0$ such that for all $\mu \in A$ it holds that $\intabx M{\infty}{x^2}{\mu(x)} < \eps$.

\end{definition}

\begin{lemma} (Relation between narrow topology and Wasserstein-$2$ metric). Proposition~7.1.5 in \cite{AmbrosioGigliSavare08} \label{lem:AGS:P715}
\[ \Wdist{\mu_n}{\mu} \rightarrow 0 \lrar \left\{ \begin{array}{l}
  \displaystyle \mu_n \weakto \mu,  \\
  \displaystyle (\mu_n) \textup{ has uniformly integrable second moments.} \\
\end{array} \right. \]

\end{lemma}


\begin{lemma} (Sufficient condition uniformly integrable second moments). \label{prop_uniformly integrable second moments_suff_cond}
Let $\mathcal A \subset \PP$ and $\varphi:[0,\infty)\to\mathbb R$ non-decreasing with $\lim_{x\to\infty}x^{-2}\varphi(x)=\infty$. If
\[
\exists \; C > 0 \; \forall \; \mu \in \mathcal A : \intabx{}{}{\varphi}\mu \leq C,
\]
then $\mathcal A$ has uniformly integrable second moments.

\end{lemma}

Now, we establish the limsup inequality. By a diagonal argument (see Lemma~3.4 in \cite{VanMeursMunteanPeletier14} for details), it is sufficient to show that for some $\mathcal Y \subset \Ptwo$ we have 
\begin{subequations}
\label{for:app:limp:cond}
\begin{alignat}2
\label{for:app:limp:cond1}
\forall \: \mu \in \mathcal Y \: \: \exists \: \mu_n \xrightarrow{\Wdistnobraces} \mu : \quad  
& \limsup_{n\to\infty}  E_n(\mu_n) \: &\leq E(\mu), \\
\forall \: \mu \in \Ptwo \: \: \exists \: \mathcal Y \ni \nu_n \xrightarrow{\Wdistnobraces} \mu : \quad  
& \limsup_{n\to\infty} E(\nu_n) &\leq E(\mu),
\label{for:app:limp:cond2}
\end{alignat}
\end{subequations}
Take $\mathcal Y$ as in \cite{GeersPeerlingsPeletierScardia13} and $\mu \in \mathcal Y$ arbitrary. For the argument below, we only require that the elements of $\mathcal Y$ are absolutely continuous, i.e.~$\rho := \ga \mu{\mathcal L} \in L^1(0, \infty)$. Due to the same arguments as in \cite{GeersPeerlingsPeletierScardia13}, the sequence $\mu_n$ defined by \eqref{eqn:rec:seq} satisfies the limsup inequality \eqref{for:app:limp:cond1}. To show that $\mu_n$ converges in $\Wdistnobraces$, we take $\varphi$ as in Lemma \ref{prop_uniformly integrable second moments_suff_cond} such that $\intabx{}{}{\varphi}\mu \leq C$. Then we put $x_{n+1}^n := \infty$ for convenience and observe
\begin{align} \notag
\frac{n}{n + 1} \intabx 0\infty{\varphi}{\mu_n}
&= \frac1{n+1} \sum_{i=0}^n \varphi \lrhaa{ x_i^n }
= \sum_{i=0}^n \varphi \lrhaa{ x_i^n } \int_{x_i^{n}}^{x_{i+1}^{n}} \rho \\ \label{for:limp:expl:recseq}
&\leq \sum_{i=0}^n \int_{x_i^{n}}^{x_{i+1}^{n}} {\varphi \rho}
= \int_0^\infty \varphi \rho.
\end{align}
Together with Lemma \ref{prop_uniformly integrable second moments_suff_cond} and Lemma \ref{lem:AGS:P715}, this result implies convergence of $\mu_n$ with respect to $\Wdistnobraces$. This completes the proof of \eqref{for:app:limp:cond1}.

To establish \ref{for:app:limp:cond2}, we take the same sequence $(\nu_n)$ as in \cite{GeersPeerlingsPeletierScardia13}. With the same strategy as before (see \eqref{for:limp:expl:recseq}) it is possible to show that $\nu_n$ converges with respect to $\Wdistnobraces$.


\providecommand{\bysame}{\leavevmode\hbox to3em{\hrulefill}\thinspace}
\providecommand{\MR}{\relax\ifhmode\unskip\space\fi MR }
\providecommand{\MRhref}[2]{%
  \href{http://www.ams.org/mathscinet-getitem?mr=#1}{#2}
}
\providecommand{\href}[2]{#2}

\end{document}